\documentclass[10pt]{amsart}
\usepackage{amsmath, amssymb, latexsym, enumerate}
\input{diagrams}

\newtheorem{thm}[equation]{Theorem}
\newtheorem{pro}[equation]{Proposition}
\newtheorem{cor}[equation]{Corollary}
\newtheorem{lem}[equation]{Lemma}

\theoremstyle{definition}
\newtheorem{exa}[equation]{Example}
\newcommand{\spandsp}{\qquad\text{and}\qquad}
\newcommand{\co}{\colon}
\newcommand{\coo}{\,\colon\,}
\newcommand{\rta}{\rightarrow}

\renewcommand{\=}{\,=\,}
\newcommand{\+}{\,+\,}

\newcommand{\isom}{\cong}
\newcommand{\sub}[1]{_{\text{${\scriptscriptstyle #1}$}}}

\newcommand{\diff}{\backslash}

\newcommand{\as}{\ensuremath{\text{$A\sub S$}}}
\newcommand{\at}{\ensuremath{\text{$A\sub T$}}}

\newcommand{\bs}{\ensuremath{\text{$B\sub S$}}}
\newcommand{\bt}{\ensuremath{\text{$B\sub T$}}}

\newcommand{\eu}{\ensuremath{\text{$E\sub U$}}}
\newcommand{\ev}{\ensuremath{\text{$E\sub V$}}}

\newcommand{\us}{\ensuremath{\text{$U\sub S$}}}
\newcommand{\ut}{\ensuremath{\text{$U\sub T$}}}
\newcommand{\vs}{\ensuremath{\text{$V\sub S$}}}
\newcommand{\vt}{\ensuremath{\text{$V\sub T$}}}

\newcommand{\rkm}[1]{\ensuremath{\rho_{\scriptscriptstyle M}(#1)}}
\newcommand{\rkmd}[1]{\ensuremath{\rho_{\scriptscriptstyle M^*}(#1)}}
\newcommand{\rkn}[1]{\ensuremath{\rho_{\scriptscriptstyle N}(#1)}}
\newcommand{\rkl}[1]{\ensuremath{\rho_{\scriptscriptstyle L}(#1)}}

\newcommand{\rkld}[1]{\ensuremath{\rho_{\scriptscriptstyle L^*}(#1)}}

\newcommand{\nlmd}[1]{\ensuremath{\nu_{\scriptscriptstyle M^*}(#1)}}
\newcommand{\nln}[1]{\ensuremath{\nu_{\scriptscriptstyle N}(#1)}}
\newcommand{\rlm}[1]{\ensuremath{\lambda_{\scriptscriptstyle M}(#1)}}

\newcommand{\frp}{\mathbin{\Box}}
\newcommand{\ci}{\ensuremath{\mathcal I}}
\newcommand{\cm}{\ensuremath{\mathcal M}}
\newcommand{\dual}{\ensuremath{^*}}

\newcommand{\ist}[1]{\ensuremath{\text{\rm Isth}(#1)}}
\newcommand{\loo}[1]{\ensuremath{\text{\rm Loop}(#1)}}

\begin{document}
\title[The Free product of Matroids]
{The Free product of Matroids}
\author{Henry Crapo}
\author{William Schmitt}
\thanks{Schmitt partially supported by NSA grant 02G-134}
\email{crapo@ehess.fr and wschmitt@gwu.edu}
\keywords{Matroid, free product, Welsh conjecture}
\subjclass[2000]{05B35, 16W30, 05A15}

\begin{abstract}
  We introduce a noncommutative binary operation on matroids, called
free product.  We show that this operation respects matroid duality, and
has the property that, given only the cardinalities, an ordered
pair of matroids may be recovered, up to isomorphism, from
its free product.  We use these results to give a short proof of Welsh's
1969 conjecture, which provides a progressive lower bound for
the number of isomorphism classes of matroids on an $n$-element set.
\end{abstract}

\maketitle

In the {\it minor coalgebra\/} of matroids
(\cite{sc:iha}, \cite{crsc:fsa}), the coproduct of a matroid $M(S)$ is
given by $\sum\sub{A\subseteq S}M|A\otimes M/A$, where $M|A$ denotes
the restriction of $M$ to $A$ and $M/A$ denotes the matroid on the set
difference $S\backslash A$ obtained by contracting $A$ from $M$.  The
product of matroids $M$ and $N$ in the dual algebra is thus a linear
combination $\sum\sub L \alpha\sub L L$ of those matroids $L$ 
having some restriction isomorphic to $M$, with
complementary contraction isomorphic to $N$.  The coefficient
$\alpha\sub L$ of $L=L(U)$ is the number of subsets $A\subseteq U$ such
that $L|A\isom M$ and $L/A\isom N$.
If the matroids having nonzero coefficient in the product
of $M$ and $N$ are ordered in the weak-map order, there is a final term
equal to a scalar multiple of the direct product $M\oplus N$, and an initial 
term equal to a scalar multiple of a matroid 
that we have elected to call the {\it free product\/} of $M$ and $N$.  

In the present short article we give an intrinsic definition of
the free product of matroids, and prove the crucial
result that, given only their cardinalities, the two
factors themselves, and even the order of the factors,
can be recovered, up to isomorphism, from the free product.
This is in sharp contrast to the behavior of direct sums, where the
failure of unique ordered factorization gave rise to a little crisis in
matroid theory, holding up the proof of Welsh's ``self-evident'' conjecture
\cite{we:bnm} for more than three decades.  He conjectured that if
there are $f_n$ isomorphism classes of matroids on an $n$-element 
set, then $f_n\cdot f_m\leq f_{n+m}$, for all $n,m\geq 0$.  
Where direct sum fails, free product succeeds; we prove the
conjecture here.

In future work we shall investigate in detail the combinatorial
properties of the free product, as well as its implications
for the minor coalgebra of matroids.

We denote the rank and nullity functions of a matroid $M(S)$ by
$\rho\sub M$ and $\nu\sub M$, respectively, and denote by
$\lambda\sub M$ the {\it rank-lack} function on $M$, given by $\rlm A=
\rho (M) - \rkm A$, for all $A\subseteq S$, where $\rho (M)
=\rkm S$ is the rank of $M$. We denote the disjoint
union of sets $S$ and $T$ by $S+T$ and the intersection $S\cap T$ by
either $S\sub T$ or $T\sub S$.  We refer the reader to Oxley's book
\cite{ox:mt} for any background on matroid theory that might be
needed.

\begin{pro}
For all matroids $M(S)$ and $N(T)$, the collection
$$
\ci = \{ A\subseteq S+T\coo\text{$A\sub S$ is independent in $M$
and $\rlm\as\geq\nln\at$}\}
$$
is the family of independent subsets of a matroid $M\frp N$ on $S+T$.
\end{pro}
\begin{proof}
Suppose that $A\in\ci$ and $B\subseteq A$.  Then $\bs\subseteq\as$ and
$\bt\subseteq\at$, and so $\rlm\bs\geq\rlm\as$ and $\nln\at\geq\nln\bt$.
Since $A\in\ci$, the set \as\ is independent in $M$ and
$\rlm\as\geq\nln\at$, from which it follows that \bs\ is independent
in $M$ and $\rlm\bs\geq\nln\bt$, that is, $B\in\ci$.  

Now suppose that $A,B\in\ci$, with $|A|<|B|$.  
We consider
three cases:  First, if $\rlm\as>\nln\at$ and $|\as|<|\bs|$, then
$\rlm{\as\cup x}\geq\nln\at$, for any $x\in\bs\diff\as$, so  if we choose
any such $x$ with $\as\cup x$ independent in $M$, then
$A\cup x =(\as\cup x)\cup\at$ belongs to \ci.
Second, if $\rlm\as>\nln\at$ and $|\at|<|\bt|$, then 
$\rlm\as\geq\nln{\at\cup y}$, for any $y\in\bt\diff\at$, and so
$A\cup y=\as\cup (\at\cup y)\in\ci$.  Finally, we consider the case
in which $\rlm\as = \nln\at$. Since $A\in\ci$, the set $\as$
is independent in $M$, so that $\rkm\as=|\as|$, and hence in this case,
$\rlm\as = \rho (M)-|\as| =\nln\at$.  
Also, since $B\in\ci$, we have $\rlm\bs=
\rho (M)-|\bs|\geq\nln\bt$, and thus
\begin{align*}
\rkn\at &\=  |\at|\,-\,\nln \at\\
& \= |\at|\+|\as|\,-\,\rho (M)\\
& \,<\, |\bt|\+|\bs|\,-\,\rho (M)\\
& \,\leq\, |\bt|\, -\, \nln\bt\\
&\=\rkn\bt.
\end{align*}
Therefore we may find an element $z$ of $\bt\diff\at$
such that $\nln{\at\cup z}=\nln \at$, and so
$A\cup z=\as\cup (\at\cup z)$ belongs to \ci.  Hence $M\frp N$ is
a matroid.
\end{proof}
We refer to the matroid $M\frp N$ as the {\it free product} of $M$ and $N$.
Note that the set of bases of $M(S)\frp N(T)$ is given by
$$
\{ A\subseteq S+T\coo\text{\as\ is independent in $M$, \at\ spans
  $N$, and $\rlm\as=\nln\at$}\},
$$
and so, in particular,  $\rho (M\frp N)=\rho (M)+\rho (N)$,
for all $M$ and $N$.  
\begin{exa}
Let $S=\{a,b,c\}$ and $T=\{d,e,f,g\}$, and suppose that $M(S)$ is a
three-point line and $N$ consists of two double points, $de$
and $fg$.  The free product $M\frp N$ is shown below:

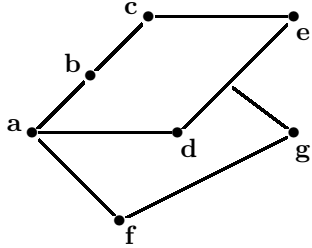
\newcommand{\eelab}[1]
{\hspace{2.ex}\text{\raisebox{.8ex}{$\bf #1$}}}
\newcommand{\elab}[1]
{\hspace{-3.1ex}\text{\raisebox{-1.5ex}{$\bf #1$}}}

\vspace{.5ex}
\begin{diagram}[PostScript=Rokicki,abut,height=1.1em,width=1.1em,tight,thick,balance]
&&&&\eelab{c}&\bullet&&\rLine&&&\bullet&\elab{\!e}\\
&&&&\ruLine(2,2)&&&&&\ruLine(4,4)\\
&&\eelab{b}&\bullet&&&&
\rdLine[leftshortfall=2.6em](4,3)\\
&&\ruLine(2,2)\\
\eelab{a}&\bullet&\rLine&&&&\bullet&\elab{d}&&&
\bullet&\elab{\!g}\\
&&\rdLine(3,3)&&&&&&&\ldLine(6,3)\\\\
&&&&\bullet&\elab{f}
\end{diagram}
\vspace{1.ex}

\noindent
The bases of $M\frp N$ are the sets of the form $A\cup B$, 
with $A\subseteq S$, $B\subseteq T$,
and either
\begin{enumerate}[(i)]
\item $A=\emptyset$\, and\, $B=T$,
\item $|A|=1$\, and\, $|B|=3$,\, or
\item $|A|=2$\, and\, $|B|=2$,\, with $B$ not equal to $\{d,e\}$ or $\{f,g\}$.
\end{enumerate}

\end{exa}
The following proposition verifies
that, in a free product $M(S)\frp N(T)$, the restriction
to $S$ and contraction by $S$ yield $M$ and $N$ as 
minors.
\begin{pro}\label{pro:unif1}
  For all matroids $M(S)$ and $N(T)$,
$$
(M\frp N)|S \= M \spandsp\ (M\frp N)/S \= N.
$$
\end{pro}
\begin{proof}
  It is immediate from the definition of independence in
 $L=M\frp N$ that $A\subseteq S$ is independent in
$L$ if and only if it is independent in $M$, and hence $L|S=M$.
Let $E\subseteq S$ be a basis of $M$.  The independent sets
of $L/S$ are those sets $B\subseteq T$ such that 
$E\cup B$ is independent in $L$.  Since $\rlm E=0$, the set
$E\cup B$ is independent in $L$ if and only if $B$ is
independent in $N$, and hence $L/S=N$.
\end{proof}

For any matroid $M=M(S)$, the rank function of the dual
matroid $M\dual$ satisfies $\rkmd B = |B|-\rho (M)+\rkm A$,
or equivalently, $\rlm A=\nlmd B$, for all $A+B=S$.
\begin{pro}\label{pro:dual}
  For all matroids $M$ and $N$, $(M\frp N)\dual = N\dual\frp M\dual$.
\end{pro}
\begin{proof}
  Suppose that $M=M(S)$, $N=N(T)$, and $A+B=S+T$, so that $A$ is
a basis for $M\frp N$ if and only if $B$ is a basis for $(M\frp N)\dual$.
Now $A$ is a basis for $M\frp N$ if and only if \as\ is independent
in $M$, \at\ spans $N$ and $\rlm\as=\nln\at$, which is true if and only
if \bs\ spans $M\dual$, \bt\ is independent in $N\dual$, and 
$\nlmd\bs=\lambda_{\scriptscriptstyle N^*}(\bt)$, that is, if
and only if $B$ is a basis for $N\dual\frp M\dual$.
\end{proof}

The next result implies that, given the size of $S$,
we can recover the rank of $M(S)$ from the free 
product $M\frp N$.
\begin{lem}\label{pro:ineq}
If $L=M(S)\frp N(T)$, then $\rho (M) \leq\rkl U$, 
for all $U\subseteq S+T$ such that $|U|=|S|$. 
\end{lem}
\begin{proof}
 Suppose that $U\subseteq S+T$, with $|U|=|S|$, and let  
$V$ be the complement of $U$ in $S+T$. Note that \vs\ is
the complement of \us\ in $S$, while \ut\ is the 
complement of \us\ in $U$, so that $|\vs|=|\ut|$.
 Let $E\subseteq\us$ be a basis for
$M|\us=L|\us$.
Since $\rlm\us\leq |\vs|=|\ut|$,
there exists $A\subseteq\ut$ with $|A|=\rlm\us$.
It follows that $E\cup A$ is independent in $L$, and thus, since
$E\cup A\subseteq U$, we have
$\rkl U\geq |E|+|A|=\rkm\us+\rlm\us=\rho (M)$.
\end{proof}
For any matroid $M$, we write \loo M\ and \ist M, respectively, for the 
sets of loops and isthmuses of $M$.  The following result shows, in
particular, that whenever $M$ is isthmusless and $N$ loopless, and we
know the size of $M$, 
then the support set of $M$, and thus the matroid $M$ itself, can
be recovered from the free product $M\frp N$
\begin{lem}\label{pro:rkineq}
  Suppose that $L=M(S)\frp N(T)$, and $U\subseteq
  S+T$ satisfies $|U|=|S|$. 
If $U$ contains a nonloop of $N$ and the
complement of $U$ in $S+T$ contains a nonisthmus of $M$, 
then $\rkl U>\rho (M)$.
\end{lem}
\begin{proof}
  Let $V$ denote the complement of $U$ in $S+T$.  Since $\vs\not\subseteq
\ist M$, and $|U|=|S|$, it follows that $\rlm\us< |\vs|=|\ut|$.
If $F$ is any subset of \ut\ containing at least one nonloop
and having size $\rlm\us +1$, then $\rlm\us\geq\nln F$, and thus
$E\cup F$ is independent in $L$, for all independent $E$ in
$M|\us$.  In particular, if $E$ is a basis for $M|\us$, then
$|E\cup F|=|E|+|F|>\rkm\us +\rlm\us = \rho (M)$,
and hence $\rkl U >\rho (M)$.
\end{proof}
Our main result, below, implies that, given the size of $M$,
we may recover $M$ and $N$ up to isomorphism from 
$M\frp N$, without conditions on $M$ and $N$.
\begin{thm}\label{thm}
If $L=M(S)\frp N(T)$, then
for any $U\subseteq S+T$ such that $|U|=|S|$ and $\rkl U=\rho (M)$,
there exist bijective weak maps $L|U\rta M$ and $L/U\rta N$.
\end{thm}
\begin{proof}
Let $V$ denote the complement of $U$ in $S+T$, 
let $f\co\vs\rta\ut$ be an arbitrary bijection, and
define $\phi\coo S+T\rta S+T$ by
$$
\phi (x)\=
\begin{cases}
  f(x), & \text{if $x\in\vs$},\\
f^{-1}(x), & \text{if $x\in\ut$},\\
x, & \text{if $x\in\us\cup\vt$}.
\end{cases}
$$
Denote by $\phi_1$ and $\phi_2$, respectively, the restrictions
$\phi|U$ and $\phi|V$, and note that $\phi_1\co U\rta S$ and
$\phi_2\co V\rta T$ are bijections.  We now show that 
$\phi_1$ and $\phi_2$ are the desired weak maps.

According to Lemma \ref{pro:rkineq}, the fact that $\rkl U=\rho
(M)$, implies that either $\vs\subseteq\ist M$ or $\ut\subseteq\loo
N$, or both; we first consider the case in which $\vs\subseteq\ist
M$.

Since $\vs\subseteq\ist M$, the bases for $M$ are the sets of
the form $A\cup\vs$, where $A$ is a basis for $M|\us$.  Now if
$A$ is independent in $M|\us$, then 
$$
\rlm A\geq |\vs|=|\ut|\geq\nln\ut,
$$
so that $A\cup\ut$ is independent in $L$, and thus the bases of
$L|U$ are the sets of the form $A\cup\ut$, where $A$ is a basis
of $M|\us$.  Hence $B=A\cup\vs$ is a basis of $M$ if and only if
$\phi^{-1}(B)=A\cup\ut$ is a basis of $L|U$, and therefore $\phi_1$
is an isomorphism from $L|U$ onto $M$.

Let $A$ be a basis for
$M|\us$, so that $B=A\cup\ut$ is a basis for $L|U$, as seen above,
and let 
$E=\eu\cup\ev$ be a basis for $N$.  In order to see that
$\phi_2\co L/U\rta N$ is a weak map, we need
to show that $\phi^{-1}(E)$ is independent in $L/U$, or equivalently,
that $\phi^{-1}(E)\cup B$ is independent in $L$.
Now $\phi^{-1}(E)=\phi^{-1}(\eu)\cup\ev$, and thus
$\phi^{-1}(E)\cup B=(\phi^{-1}(\eu)\cup A)\cup (\ev\cup\ut)$
is independent in $L$ if and only if
$\nln{\ev\cup\ut}\leq\rlm{\phi^{-1}(\eu)\cup A}$.
Since $\phi^{-1}(\eu)\subseteq\vs\subseteq\ist M$, and $A$ is
a basis for $M|\us$, we have
\begin{align*}
  \rlm{\phi^{-1}(\eu)\cup A} &\=|\vs|-|\phi^{-1}(\eu)|\\
&\=|\vs|-|\eu|.
\end{align*}
On the other hand, since  $E\subseteq\ev\cup\ut$, and $E$ is a basis 
for $N$, we have
\begin{align*}
  \nln{\ev\cup\ut}&\= |\ev\cup\ut|-\rkn{\ev\cup\ut}\\
&\= |\ev|+|\ut|-|E|\\
&\= |\ut|-|\eu|.
\end{align*}
Hence $\nln{\ev\cup\ut}=\rlm{\phi^{-1}(\eu)\cup A}$, and so
$\phi^{-1}(E)\cup B$ is independent in (and in fact is a basis
for) $L$.  Thus, in the case that $\vs\subseteq\ist M$, we have
that $\phi_1\co L|U\rta M$ and $\phi_2\co L/U\rta N$ are
weak maps.

Now suppose that $\ut\subseteq\loo N=\ist{N\dual}$. 
By Proposition \ref{pro:dual}, we have 
$L\dual = N\dual\frp M\dual$, and since
$|U|=|S|$ and $\rkl U=\rho (M)$, it follows that $|V|=
|T|$ and $\rkld V=\rho (N\dual)$.  Interchanging the roles
of $M$, $N$, and $U$, respectively, with those of $N\dual$,
$M\dual$ and $V$ in the above, we obtain that 
$\phi_2\co L\dual|V=(L/U)\dual\rta N\dual$ and $\phi_1\co
L\dual/V=(L|U)\dual\rta M\dual$ are weak maps, and since
$\rho ((L/U)\dual)=\rho (N\dual)$ and $\rho ((L|U)\dual)=
\rho (M\dual)$, this implies
that $\phi_1\co L|U\rta M$ and $\phi_2\co L/U\rta N$
are weak maps.  
\end{proof}

\begin{cor}\label{big}
  If $M(S)\frp N(T)\isom P(U)\frp Q(V)$, where $|S|=|U|$, 
then $M\isom P$ and $N\isom Q$.
\end{cor}
\begin{proof}
  Choosing an isomorphism from $M\frp N$ to $P\frp Q$, and
relabelling if necessary, we can assume that $M\frp N=P\frp Q$.
Since $\rho (M)=\rho\sub{P\frp Q}(S)$ and $\rho (P)=\rho
\sub{M\frp N}(U)$, it follows from Lemma \ref{pro:ineq}
that $\rho (M)=\rho (P)$.  We may thus apply Theorem \ref{thm}
to obtain bijective weak maps $M\rta P$, $N\rta Q$, $P\rta M$ and
$Q\rta N$; hence $M\isom P$ and $N\isom Q$.
\end{proof}
 
We thus obtain the following result, which was conjectured by Welsh in 
\cite{we:bnm}:
\begin{cor}
  If $f_n$ denotes the number of nonisomorphic matroids on a
set of size $n$, then $f_{n+m}\geq f_n\cdot f_m$, for all $n,m\geq 0$.
\end{cor}
\begin{proof}
Denote by $\cm (n)$ be the set of all isomorphism
classes of matroids on a set of $n$ elements, so that 
$|\cm (n)|=f_n$, for all $n\geq 0$, and write
$[M]$ for the isomorphism class of a matroid $M$.  Corollary
\ref{big} says precisely that, for all $n,m\geq 0$, 
the map $\cm (n)\times\cm (m)\rta\cm (n+m)$ given by
$([M],[N])\mapsto [M\frp N]$ is injective, from which
the inequality follows.
\end{proof}

\end{document}